\documentclass{amsart}

\usepackage{fancyhdr}
\usepackage{lastpage}
\usepackage{stmaryrd,yhmath}

\newcommand{\bdis}{\begin{displaymath}}
\newcommand{\edis}{\end{displaymath}}
\newcommand{\be}{\begin{equation}}
\newcommand{\ee}{\end{equation}}
\newcommand{\mbb}{\mathbb}
\newcommand{\mcal}{\mathcal}

\newcommand{\vp}{\varphi}

\newcommand{\zf}{\zeta\left(\frac{1}{2}+it\right)}

\theoremstyle{definition}

\theoremstyle{remark}
\newtheorem{remark}[]{Remark}

\newtheorem*{mydef1}{{\bf Theorem}}

\newtheorem*{mydef4}{{\bf Corollary}}

\numberwithin{equation}{section}

\begin{document}

\title{Jacob's ladders and the oscillations of the function $|\zeta(1/2+it)|^2$ around its mean-value; law of the almost exact equality of
corresponding areas}

\author{Jan Moser}

\address{Department of Mathematical Analysis and Numerical Mathematics, Comenius University, Mlynska Dolina M105, 842 48 Bratislava, SLOVAKIA}

\email{jan.mozer@fmph.uniba.sk}


\begin{abstract}
The oscillations of the function $Z^2(t),\ t\in [0,T]$ around the main part $\sigma(T)$ of its mean-value are studied in this paper. It is proved that
an almost equality of the corresponding areas holds true. This result cannot be obtained by methods of Balasubramanian, Heath-Brown and Ivic.
\end{abstract}

\maketitle

\section{Introduction}

\subsection{}

The Titchmarsh-Kober-Atkinson (TKA) formula
\be \label{1.1}
\int_0^\infty Z^2(t)e^{-2\delta t} {\rm d}t=\frac{c-\ln (4\pi\delta)}{2\sin\delta}+\sum_{n=0}^N c_n\delta^n+\mcal{O}(\delta^{N+1})
\ee
(see \cite{17}, p. 131) remained as an isolated result for the period of 56 years. We have discovered (see \cite{5}) the nonlinear integral equation
\be \label{1.2}
\int_0^{\mu[x(T)]}Z^2(t)e^{-\frac{2}{x(T)}t}{\rm d}t=\int_0^T Z^2(t){\rm d}t
\ee
in which the essence of the TKA formula is encoded. Namely, we have shown in \cite{5} that the following almost exact expression for the Hardy-Littlewood
integral
\be \label{1.3}
\int_0^T Z^2(t){\rm d}t=\frac{\vp(T)}{2}\ln\frac{\vp(T)}{2}+(c-\ln (2\pi))\frac{\vp(T)}{2}+c_0+\mcal{O}\left(\frac{\ln T}{T}\right)
\ee
takes place, where $\vp(T)$ is the Jacob's ladder, i.e. an arbitrary solution to the nonlinear integral equation (\ref{1.2}).

\begin{remark}
Our formula (\ref{1.3}) for the Hardy-Littlewood integral
\be \label{1.4}
\int_1^T\left|\zf\right|^2{\rm d}t=\int_1^T Z^2(t){\rm d}t
\ee
has been obtained after the time period of 90 years since this integral appeared in 1918 with the first result
\bdis
\int_1^T\left|\zf\right|^2{\rm d}t\sim T\ln T
\edis
(see \cite{3}, pp. 122, 151-156).
\end{remark}

\subsection{}

Let us remind that

\begin{itemize}
\item[(A)] The Good's $\Omega$ - theorem (see \cite{2}) implies for the Balasubramanian formula (see \cite{1})
\be \label{1.5}
\int_0^T Z^2(t){\rm d}t\sim T\ln T+ (2c-1-\ln 2\pi)T+R(T),\ R(T)=\mcal{O}(T^{1/3+\epsilon})
\ee
that
\be \label{1.6}
\limsup_{T\to\infty}|R(T)|=+\infty ,
\ee
i.e. the error term in (\ref{1.5}) is unbounded at $T\to\infty$.

\item[(B)] In the case of our formula (\ref{1.3}) the error term definitely tends to zero
\be \label{1.7}
\lim_{T\to\infty} r(T)=0; \qquad r(T)=\mcal{O}\left(\frac{\ln T}{T}\right) ,
\ee
i.e. our formula is almost exact (see \cite{5}).

\end{itemize}

\begin{remark}
In this paper the geometric interpretation of (\ref{1.6}) and (\ref{1.7}) is obtained.
\end{remark}

\subsection{}

For the mean-value of the function $\left|\zf\right|^2=Z^2(t)$, where
\bdis
Z(t)=e^{i\vartheta(t)}\zf,\quad \vartheta(t)=-\frac t2\ln\pi+\text{Im}\ln\Gamma\left(\frac 14+i\frac t2\right),
\edis
we obtain from (\ref{1.3})
\be \label{1.8}
\frac 1T\int_0^T Z^2(t){\rm d}t=\frac{\vp(T)}{2T}\ln\frac{\vp(T)}{2}+(c-\ln 2\pi)\frac{\vp(T)}{2T}+\frac{c_0}{T}+\mcal{O}\left(\frac{\ln T}{T}\right) .
\ee
Let
\be \label{1.9}
\sigma(T)=\frac{\vp(T)}{2T}\ln\frac{\vp(T)}{2}+(c-\ln 2\pi)\frac{\vp(T)}{2T}+\frac{c_0}{T}
\ee
denote the main part of the mean-value (\ref{1.8}). In this paper the oscillation of the values of the function $Z^2(t),\ t\in [0,T]$ around
the main part $\sigma(T)$ of its mean-value are studied.

\begin{remark}
The main result of this paper is the following statement: the areas of the figures corresponding to the parts of the graph of the function
$Z^2(t),\ t\in [0,T]$ given by inequalities $Z^2(t)\geq \sigma(T)$ and $Z^2(t)\leq \sigma(T)$, respectively, are almost exactly equal.
\end{remark}

This paper is a continuation of the series \cite{5}-\cite{16}.

\section{Result}

\subsection{}

Let (see (\ref{1.9}))

\be \label{2.1}
\begin{split}
S^+(T) & =  \{ t:\ Z^2(t)\geq \sigma(T),\ t\in [0,T]\} , \\
S^-(T) & =  \{ t:\ Z^2(t)< \sigma(T),\ t\in [0,T]\}
\end{split}
\ee
and
\be \label{2.2}
\begin{split}
\Pi^+(T) & = \{ (t,y):\ \sigma(T)\leq y\leq Z^2(t),\ t\in S^+(T)\} , \\
\Pi^-(T) & = \{ (t,y):\ Z^2(T)\leq y\leq \sigma(t),\ t\in S^-(T)\} ,
\end{split}
\ee
i.e. $\Pi^+$ is the figure that corresponds to the parts of the graph of $y=Z^2(t),\ t\in [0,T]$ lying above the segment $y=\sigma(T)$ and
similarly $\Pi^-$ corresponds to the parts of the graph lying under that segment. Let $m\{\Pi^+(T)\},\ m\{\Pi^-(T)\}$ denote
measures of corresponding figures, i.e.

\be \label{2.3}
\begin{split}
m\{\Pi^+(T)\} & = \int_{S^+(T)} \{ Z^2(t)-\sigma(T)\}{\rm d}t , \\
m\{\Pi^-(T)\} & = \int_{S^-(T)} \{ \sigma(T)-Z^2(t)\}{\rm d}t .
\end{split}
\ee
The following theorem holds true.

\begin{mydef1}
First of all, we have the formula
\be \label{2.4}
m\{\Pi^+(T)\}=m\{\Pi^-(T)\}+\mcal{O}\left(\frac{\ln T}{T}\right)
\ee
(see (\ref{1.3}), (\ref{1.9}), (\ref{2.1})-(\ref{2.3})). Next, the structure of the formula (\ref{2.4}) is as follows: there are the functions
$\eta_1(T),\ \eta_2(T)$ that the following formulae
\be \label{2.5}
\begin{split}
m\{\Pi^+(T)\} & = \frac{1+o(1)}{2\pi^2}\frac{T\ln^4 T}{\eta_1-\eta_2}-\frac{\eta_2}{\eta_1-\eta_2}\mcal{O}\left(\frac{\ln T}{T}\right), \\
m\{\Pi^-(T)\} & = \frac{1+o(1)}{2\pi^2}\frac{T\ln^4 T}{\eta_1-\eta_2}-\frac{\eta_1}{\eta_1-\eta_2}\mcal{O}\left(\frac{\ln T}{T}\right)
\end{split}
\ee
hold true, and
\be \label{2.6}
AT^{2/3}\ln^4 T< m\{\Pi^+(T)\},m\{\Pi^-(T)\}<AT\ln T .
\ee
In addition to (\ref{2.6}): on the Lindel\" of hypothesis
\be \label{2.7}
m\{\Pi^+(T)\},m\{\Pi^-(T)\}>A(\epsilon)T^{1-\epsilon} ,
\ee
and on Riemann hypothesis
\be \label{2.8}
m\{\Pi^+(T)\},m\{\Pi^-(T)\}>T^{1-\frac{A}{\ln\ln T}} .
\ee
\end{mydef1}

\begin{mydef4}
We have by (\ref{2.5}), (\ref{2.6})
\be \label{2.9}
\eta_1(T)-\eta_2(T)>A \ln^3 T .
\ee
\end{mydef4}

\begin{remark}
Since from (\ref{2.4})
\be \label{2.10}
\lim_{T\to\infty}[m\{\Pi^+(T)\}-m\{\Pi^-(T)\}]=0
\ee
follows then we have the almost exact equality of the areas $m\{\Pi^+(T)\}$ and $m\{\Pi^-(T)\}$.
\end{remark}

\subsection{}

In the case of the Balasubramanian formula (\ref{1.5}) we have (comp. (\ref{1.3}), (\ref{1.9}))
\bdis
\sigma_1(T)=\ln T+2c-1-\ln 2\pi .
\edis
Let $S^+_1(T), S^-_1(T), \Pi^+_1(T), \Pi^-_1(T), m\{\Pi^+_1(T)\}, m\{\Pi^-_1(T)\}$ correspond to $\sigma_1(T)$ similarly to (\ref{2.1})-(\ref{2.3}).
Then from (\ref{1.5}) we obtain
\be \label{2.11}
m\{\Pi^+_1(T)\}=m\{\Pi^-_1(T)\}+\mcal{O}(T^{1/3+\epsilon}),\ T\to\infty ,
\ee
and (see (\ref{1.6})
\be \label{2.12}
\limsup_{T\to\infty}|m\{\Pi^+_1(T)\}-m\{\Pi^-_1(T)\}|=+\infty .
\ee

\begin{remark}
The following holds true:
\begin{itemize}
\item[(A)] Our formula (\ref{1.3}) which has been obtained by means of the Jacob's ladders leads to the almost exact equality of the areas
(see (\ref{2.4}), (\ref{2.10}).
\item[(B)] The Balasubramanian formula (\ref{1.5}) which has been obtained by means of estimation of trigonometric sums leads to the formula
(\ref{2.11}) that possesses quite large uncertainty (\ref{2.11}) and this error term cannot be removed.
\end{itemize}
\end{remark}

\section{Proof of Theorem}

\subsection{}

We obtain from (\ref{1.3}), (\ref{2.1})
\be \label{3.1}
\int_{S^+(T)}\{ Z^2(t)-\sigma(T)\}{\rm d}t+\int_{S^-(T)}\{ Z^2(t)-\sigma(T)\}{\rm d}t=\mcal{O}\left(\frac{\ln T}{T}\right)
\ee
and from (\ref{3.1}) by (\ref{2.3}) the formula
\be \label{3.2}
m\{\Pi^+(T)\}-m\{\Pi^-(T)\}=\mcal{O}\left(\frac{\ln T}{T}\right)
\ee
follows, i.e. (\ref{2.4}).

\subsection{}

Next, from the Ingham formula (see \cite{4}, p. 277, \cite{17}, p. 125)
\be \label{3.3}
\int_0^T Z^4(t){\rm d}t=\frac{1}{2\pi^2}T\ln^4T+\mcal{O}(T\ln^3T)
\ee
we obtain (see (\ref{1.9})
\be \label{3.4}
\int_0^T\{ Z^4(t)-\sigma^2(T)\}{\rm d}t=\frac{1}{2\pi^2}T\ln^4T-T\sigma^2(T)+\mcal{O}(T\ln^3T) .
\ee
Since ($\vp(T)\sim T$)
\bdis
T\sigma^2(T)=\mcal{O}\left\{\frac{\vp^2(T)}{T}\ln^2\frac{\vp(T)}{2}\right\}=\mcal{O}(T\ln^2T) ,
\edis
then from (\ref{3.4}) the formula
\be \label{3.5}
\int_0^T\{ Z^4(t)-\sigma^2(T)\}{\rm d}t=\frac{1+o(1)}{2\pi^2}T\ln^4T
\ee
follows.

\subsection{}

Since $Z^4(t)-\sigma^2(T)=(Z^2-\sigma)(Z^2+\sigma)$ and $Z^2(t)-\sigma(T)$ is always of the same sign on $S^+(T)$ and on $S^-(T)$, respectively, then
from (\ref{3.5}) we obtain (see (\ref{2.3}))
\be \label{3.6}
\eta_1(T)m\{\Pi^+(T)\}-\eta_2(T)m\{\Pi^-(T)\}=\frac{1+o(1)}{2\pi^2}T\ln^4T ,
\ee
where $\eta_1=\eta_1(T),\ \eta_2=\eta_2(T)$ are the mean-values of $Z^2(t)+\sigma(T)$ relatively to the values of the functions
$Z^2(t)-\sigma(T)$ and $\sigma(T)-Z^2(t)$, respectively on the sets $S^+(T)$ and $S^-(T)$, respectively. It is clear that
\be \label{3.7}
A\ln T <\eta_1(T),\eta_2(T)<A T^{1/3}
\ee
(see (\ref{1.9}); $|Z(t)|<t^{1/6}$). Next, $\eta_1(T)\not=\eta_2(T)$ is also true. Since if $\eta_1=\eta_2$ then by (\ref{3.2}), (\ref{3.6}), (\ref{3.7})
we would have the contradiction. Hence, from the simple system of linear equations (\ref{3.2}), (\ref{3.6}) we obtain (\ref{2.5}).

\subsection{}

Since
\be \label{3.8}
0<\eta_1-\eta_2<\eta_1<AT^{1/3}
\ee
(see (\ref{2.5}), (\ref{3.7})) then we obtain from (\ref{2.5}) the lower estimates in (\ref{2.6}). Next we have (see (\ref{1.9}), (\ref{2.3})
\bdis
m\{\Pi^+(T)\},m\{\Pi^-(T)\}<\int_0^T\{ Z^2(t)-\sigma(T)\}{\rm d}t<AT\ln T
\edis
i.e. the upper estimates in (\ref{2.6}) hold true.

\subsection{}

Following the Lindel\" of and the Riemann conjectures the estimates
\bdis
Z^2(t)< A(\epsilon)t^\epsilon,\quad t^{\frac{A}{\ln\ln t}}
\edis
take place correspondingly and then the conditional estimates (\ref{2.7}) and (\ref{2.8}) follow.

\thanks{I would like to thank Michal Demetrian for helping me with the electronic version of this work.}

\end{document}